\documentclass[12pt,oneside,reqno]{amsart}
\usepackage{graphicx}
\usepackage{mathrsfs}
\usepackage{stmaryrd}
\usepackage{amsfonts}
\usepackage{enumerate,amsmath,amssymb,amsthm}
\newcommand{\dif}{\mathrm{d}}
\newcommand{\me}{\mathrm{e}}
\newcommand{\be}{\begin{eqnarray}}
\newcommand{\ee}{\end{eqnarray}}
\newcommand{\ce}{\begin{eqnarray*}}
\newcommand{\de}{\end{eqnarray*}}
\newtheorem{theorem}{Theorem}[section]
\newtheorem{lemma}[theorem]{Lemma}
\newtheorem{remark}[theorem]{Remark}
\newtheorem{definition}[theorem]{Definition}
\newtheorem{proposition}[theorem]{Proposition}
\newtheorem{Example}[theorem]{Example}
\newtheorem{corollary}[theorem]{Corollary}
\def\e{\varepsilon}
\def\t{\theta}

\def\[{{\Big[}}
\def\]{{\Big]}}
\def\<{{\langle}}
\def\>{{\rangle}}
\def\({{\Big(}}
\def\){{\Big)}}

\def\no{\nonumber}
\def\bt{\begin{theorem}}
\def\et{\end{theorem}}
\def\bl{\begin{lemma}}
\def\el{\end{lemma}}
\def\br{\begin{remark}}
\def\er{\end{remark}}
\def\bx{\begin{Example}}
\def\ex{\end{Example}}
\def\bd{\begin{definition}}
\def\ed{\end{definition}}
\def\bp{\begin{proposition}}
\def\ep{\end{proposition}}
\def\bc{\begin{corollary}}
\def\ec{\end{corollary}}

\def\cB{{\mathcal B}}

\def\cL{{\mathcal L}}
\def\cM{{\mathcal M}}

\def\cP{{\mathcal P}}

\def\mE{{\mathbb E}}

\def\mL{{\mathbb L}}

\def\mN{{\mathbb N}}

\def\mP{{\mathbb P}}
\def\mQ{{\mathbb Q}}
\def\mR{{\mathbb R}}
\def\mS{{\mathbb S}}

\def\mU{{\mathbb U}}

\def\sA{{\mathscr A}}
\def\sB{{\mathscr B}}

\def\sF{{\mathscr F}}

\def\sL{{\mathscr L}}

\def\sN{{\mathscr N}}

\def\geq{\geqslant}
\def\leq{\leqslant}

\begin{document}

\allowdisplaybreaks

\title{Limit theorems of SDEs driven by L\'evy processes and application to nonlinear filtering problems*}

\author{Huijie Qiao}

\thanks{{\it AMS Subject Classification(2010):} 60H10, 60G51}

\thanks{{\it Keywords:} L\'evy processes, non-local Fokker-Planck equations, superposition principles, nonlinear filtering problems, the robustness}

\thanks{*This work was partly supported by NSF of China (No. 11001051, 11371352) and China Scholarship Council under Grant No. 201906095034.}

\subjclass{}

\date{}

\dedicatory{School of Mathematics,
Southeast University\\
Nanjing, Jiangsu 211189,  China\\
Department of Mathematics, University of Illinois at
Urbana-Champaign\\
Urbana, IL 61801, USA\\
hjqiaogean@seu.edu.cn}

\begin{abstract} 
In this paper we study
the convergence of solutions for (possibly degenerate) stochastic differential equations driven by L\'evy processes, when the coefficients converge in some appropriate sense. 
First, we prove, by means of a superposition principle, a limit theorem of stochastic differential equations 
driven by L\'evy processes. 
Then we apply the result to a type of nonlinear filtering problems and obtain the convergence of the nonlinear filterings.
\end{abstract}

\maketitle \rm

\section{Introduction}

Fix $T>0$ and consider the following stochastic differential equation (SDE in short) driven by a L\'evy process on $\mR^d$:
\be
\dif X_t=b(t,X_t)\dif t+\sigma(t,X_t)\dif B_t+f(t,X_{t-})\dif L_t, \quad t\in[0,T],
\label{Eq1}
\ee
where $(B_t)_{t\in[0,T]}$ is an $m$-dimensional Brownian motion and 
$(L_t)_{t\in[0,T]}$ is a $d$-dimensional pure jump L\'evy process 
with the L\'evy measure $\nu_1$. 
The coefficients $b: [0,T]\times\mR^d\mapsto\mR^d$, $\sigma: [0,T]\times\mR^d\mapsto\mR^d\times\mR^m$, and $f: [0,T]\times\mR^d\mapsto\mR$ are Borel measurable. 
Up to now, there have been many papers dealing with Eq.(\ref{Eq1}). We mention some
of these below.
In \cite{da}, Applebaum introduced some general theory, such as well-posedness and stochastic flows under Lipschitz conditions. 
Jacod collected a lot of results about the martingale problems in \cite{jj}. 
Later, Jacod and Shiryaev \cite{jjas} studied the limit theorems of Eq.(\ref{Eq1}) under Lipschitz conditions. Recently, Qiao and Zhang \cite{q4} proved that the solutions form a homeomorphism flow under non-Lipschitz conditions. Qiao and Duan \cite{qd} investigated the nonlinear filtering problems about Eq.(\ref{Eq1}) under non-Lipschitz conditions. Very recently, R\"ockner, Xie and Zhang \cite{RXZ} combined Eq.(\ref{Eq1}) with the non-local Fokker-Planck equation (\ref{FPE1}), and proved 
a one-to-one correspondence between martingale solutions of Eq.(\ref{Eq1}) 
and weak solutions of Eq.(\ref{FPE1}). 

The first goal of this paper is to apply
the result in \cite{RXZ} to a sequence of SDEs like Eq.(\ref{Eq1}). 
More precisely, we consider
the following sequence of SDEs driven by L\'evy processes: 
\be
\dif X_t=b(t,X_t)\dif t+\sigma(t,X_t)\dif B_t+\gamma g(t, X_{t-}) \dif L_t, \qquad t\in[0,T],
\label{Eq50}
\ee
where $\gamma\in\mR$ and $g: [0,T]\times\mR^d\mapsto\mR$ is Borel measurable, and for any $n\in\mN$,
\be
\dif X_t^n=b^n(t, X_t^n)\dif t+\sigma^n(t, X_t^n)\dif B_t+\gamma^n g(t, X^n_{t-}) \dif L_t, \quad t\in[0,T],
\label{Eq60}
\ee
where $b^n: [0,T]\times\mR^d\mapsto\mR^d$, $\sigma^n: [0,T]\times\mR^d\mapsto \mR^d\times\mR^m$ are Borel measurable functions and $\{\gamma^n\}$ is a real sequence. When $b^n\rightarrow b, a^n\rightarrow a, \gamma^n\rightarrow\gamma$ in some sense, where $a^n:=\frac{1}{2}\sigma^n\sigma^{n*}$ and $a:=\frac{1}{2}\sigma\sigma^{*}$, we prove that a martingale solution of Eq.(\ref{Eq60}) weakly converges to that of Eq.(\ref{Eq50}) through 
the superposition principle in \cite{RXZ}. In this paper, $\sigma$ can be degenerate. 

Our second aim is to apply the above result to a type of nonlinear filtering problems. Let us explain this in detail. Given the filtered probability space $(\Omega, \mathscr{F}, \{\mathscr{F}_t\}_{t\in[0,T]},\mP)$. 
Let the Brownian motion $B_{\cdot}$ and the L\'evy process $L_{\cdot}$ be defined
on $(\Omega, \mathscr{F}, \{\mathscr{F}_t\}_{t\in[0,T]},\mP)$. Consider a sequence of observation processes as follows:
\ce
&&Y_t=W_t+\int_0^t h(X_s)\dif s+\int_0^t\int_{\mU_0}u\tilde{N}_{\lambda}(\dif s, \dif u)+\int_0^t\int_{\mR^k\setminus\mU_0}u N_{\lambda}(\dif s, \dif u),\\
&&Y^n_t=W_t+\int_0^t h(X^n_s)\dif s+\int_0^t\int_{\mU_0}u\tilde{N}_{\lambda}(\dif s, \dif u)+\int_0^t\int_{\mR^k\setminus\mU_0}u N^n_{\lambda}(\dif s, \dif u),
\de
where $W_{\cdot}$ is a $k$-dimensional Brownian motion and $N_{\lambda}(\dif t,\dif u), N^n_{\lambda}(\dif t,\dif u)$ are two random measures with predictable compensators $\lambda(X_t,u)\dif t\nu_2(\dif u)$ and $\lambda(X^n_t,u)\dif t\nu_2(\dif u)$, respectively. Here the function $\lambda: \mR^d\times\mR^k\mapsto(0,1)$ is Borel measurable and $\nu_2$ is a $\sigma$-finite measure defined on $\mR^k$ with $\nu_2(\mR^k\setminus\mU_0)<\infty$ and $\int_{\mU_0}|u|^2\,\nu_2(\dif u)<\infty$ for a fixed $\mU_0\in\sB(\mR^k)$. $h: \mR^d\mapsto\mR^k$ is Borel measurable. Set
\ce
\pi_t(\phi):=\mE[\phi(X_t)|\mathscr{F}_t^{Y}], \quad \pi^n_t(\phi):=\mE[\phi(X^n_t)|\mathscr{F}_t^{Y^n}], \quad \phi\in\cB(\mR^d),
\de
where $\mathscr{F}_t^{Y} \triangleq\sigma(Y_s:
 0\leq s \leq t)$ and $\mathscr{F}_t^{Y^n} \triangleq\sigma(Y^n_s:
 0\leq s \leq t)$. We show that $\pi^n$ also weakly converges to $\pi$ as $X^n$ weakly converges to $X$.
 
Here we make some comments about our results. 
 First, if we specially take $L_t=\int_0^t\int_{\mR^d}uN(\dif s \dif u)$ in Eq.(\ref{Eq50}) and Eq.(\ref{Eq60}), Theorem \ref{limit1} overlaps with \cite[Theorem 3.1]{q0}. Second, if $\gamma=\gamma^n=0$ and $Y, Y^n$ have no jump parts, Theorem \ref{filcon} is just \cite[Theorem 9.4 (b)]{bkk1} and \cite[Theorem 3.3 (b)]{bkk2}. Therefore, our results are more general.

The content is arranged as follows. In the next section, we define martingale solutions for SDEs driven by L\'evy processes and weak solutions of the Fokker-Planck equations(FPEs in short). The superposition principle for SDEs driven by L\'evy processes and non-local FPEs and the stochastic Gronwall inequality are also introduced in the section. We state and prove a limit theorem in Section \ref{add}. In Section \ref{robunon}, the nonlinear filtering problems are introduced and then the convergence of nonlinear filterings is proved. Finally, we show Remark \ref{equi} in the appendix.

The following convention will be used throughout the paper: $C$ with
or without indices will denote different positive constants whose values may change from one place to another.

\section{Preliminary}\label{pre}

\subsection{Notation}

In this subsection, we introduce some  notation used in the sequel. 

We use $\mid\cdot\mid$ and $\parallel\cdot\parallel$  for the norms of vectors and matrices, respectively. We use $\langle\cdot$ , $\cdot\rangle$ to denote the scalar product in $\mR^d$. 

Let $\cB(\mR^d)$ denote the set of all real-valued uniformly bounded $\mathscr{B}(\mR^d)$-measurable functions on $\mR^d$. $C^2(\mR^d)$ stands for the space of continuous functions on $\mR^d$ which have continuous partial derivatives of order up to $2$, and $C_b^2(\mR^d)$ stands for the subspace of $C^2(\mR^d)$, consisting of functions whose derivatives up to order 2 are bounded. $C_c^2(\mR^d)$ is the collection of all functions in $C^2(\mR^d)$ with compact support and $C_c^\infty(\mR^n)$ denotes the collection of all real-valued $C^\infty$ functions of compact
support.

Let $\cP({\mR^d})$ be the space of all probability measures on $\sB(\mR^d)$, equipped with the topology of weak convergence. 

\subsection{Martingale solutions for SDEs driven by L\'evy processes}

In this subsection, we define martingale solutions for SDEs driven by L\'evy processes. 

By the L\'evy -It\^o theorem (\cite{sa}), we know that Eq.(\ref{Eq1}) can be rewritten as
\ce
\dif X_t&=&b(t,X_t)\dif t+\sigma(t,X_t)\dif B_t+\int_{|f(t,X_{t-})z|\leq l}f(t,X_{t-})z\tilde{N}(\dif t,\dif z)\\
&&+\int_{|f(t,X_{t-})z|>l}f(t,X_{t-})zN(\dif t,\dif z),
\de
where $l\geq 0$ is a constant, $N(\dif t,\dif z)$ is the Poission random measure associated with $(L_t)_{t\in[0,T]}$ and $\tilde{N}(\dif t,\dif z):=N(\dif t,\dif z)-\nu_1(\dif z)\dif t$. Moreover, the infinitesimal generator of $X_{\cdot}$ is formally expressed as
\ce
(\sL_t\phi)(x)&:=&a_{ij}(t,x)\partial_{ij}\phi(x)+b_i(t,x)\partial_i\phi(x)+\int_{\mR^d}\[\phi(x+u)-\phi(x)-I_{|u|\leq l}u^i\partial_i\phi(x)\]\nu^f_{t,x}(\dif u)\\
&=:&(\sA_t\phi)(x)+(\sB_t\phi)(x)+(\sN^f_t\phi)(x), \quad \phi\in C_b^2(\mR^d),
\de
where $\nu^f_{t,x}(A):=\int_{\mR^d}I_A(f(t,x)z)\nu_1(\dif z)$ for any $A\in\sB(\mR^d)$. 

Besides, let $D^d_T:=D([0,T], \mR^d)$ be the set of all the c\`adl\`ag functions from $[0, T]$ to $\mR^d$. $w$ stands for a generic element in $D^d_T$. We equip $D^d_T$ with the Skorokhod topology and then $D^d_T$ is a Polish space. For any $t\in[0,T]$, set 
$$
e_t: D_T\rightarrow\mR^d, \quad e_t(w)=w_t, \quad w\in D_T.
$$
Let $\cB_t:=\sigma\{w_s: s\in[0,t]\}$,  $\bar{\cB}_t:=\cap_{s>t}\cB_s$, and $\cB:=\cB_T$. In the following, we define martingale solutions of Eq.(\ref{Eq1}).(c.f.\cite{jjas, sv})

\bd(Martingale solutions)\label{martsolu}
For $\mu_0\in\cP(\mR^d)$ and $0\leq s<T$. A probability measure $\mQ$ on $(D^d_T, \cB)$ is called a martingale solution of Eq.(\ref{Eq1}) with the initial law $\mu_0$ at time $s$, if

(i) $\mQ(w_t=w_s, t\in[0,s])=1$ and $\mQ\circ e^{-1}_s=\mu_0$,

(ii) For any $\phi\in{C_c^2(\mR^d)}$,
\be
\cM_t^\phi&:=&\phi(w_t)-\phi(w_s)-\int_s^t(\sL_r\phi)(w_r)\dif r
\label{eq2}
\ee
is a $(\bar{\cB}_t)_{t\in[0,T]}$-adapted martingale under the probability measure $\mQ$. The uniqueness of the martingale solutions to Eq.(\ref{Eq1}) means that, if $\mQ, \tilde{\mQ}$ are two martingale solutions to Eq.(\ref{Eq1}) with $\mQ\circ e_s^{-1}=\tilde{\mQ}\circ e_s^{-1}$, then $\mQ\circ e_t^{-1}=\tilde{\mQ}\circ e_t^{-1}$ for any $t\in[s,T]$.
\ed

Now, we assume:
\begin{enumerate}[(${\bf H}^1_{b,\sigma}$)]
\item 
There is a constant $C_1\geq 0$ such that for all $(t,x)\in[0,T]\times\mR^d$,
\ce
|b(t,x)|+\|\sigma(t,x)\|\leq C_1(1+|x|).
\de
\end{enumerate}
\begin{enumerate}[(${\bf H}^s_{f}$)]
\item 
For all $(t,x)\in[0,T]\times\mR^d$,
\ce
\int_{\mR^d}I_{B_l}(f(t,x)z)|f(t,x)z|^2\nu_1(\dif z)<\infty, ~\mbox{and}~\int_{\mR^d}I_{B_l}(f(t,x)z)|f(t,x)z|^2\nu_1(\dif z)\leq C_2(1+|x|^2),
\de
where $B_l:=\{y\in\mR^d; |y|\leq l\}$ and $C_2\geq 0$ is a constant independent of $t, x$.
\end{enumerate}
\begin{enumerate}[(${\bf H}^l_{f}$)]
\item 
For all $(t,x)\in[0,T]\times\mR^d$,
\ce
\int_{\mR^d}I_{B^c_l}(f(t,x)z)\nu_1(\dif z)<\infty, ~\mbox{and}~\int_{\mR^d}I_{B^c_l}(f(t,x)z)\log\left(1+\frac{|f(t,x)z|}{1+|x|}\right)\nu_1(\dif z)\leq C_3,
\de
where $B^c_l:=\{y\in\mR^d; |y|>l\}$ and $C_3\geq 0$ is a constant independent of $t, x$.
\end{enumerate}

\br\label{equi}
Under $({\bf H}^1_{b,\sigma})$, $({\bf H}^s_{f})$ and $({\bf H}^l_{f})$, it can be justified that (ii) in Definition \ref{martsolu} is equivalent to the following condition: for any $\phi\in{C^2(\mR^d)}$ with $|\phi(x)|\leq C\log(2+|x|)$, $\cM_t^\phi$ is a local $(\bar{\cB}_t)_{t\in[0,T]}$-adapted martingale under the probability measure $\mQ$. 

For the readers' convenience, we put the verification in the appendix.
\er

\br\label{rewrcond}
(i)  By $({\bf H}^1_{b,\sigma})$ and $a(t,x)=\frac{1}{2}\sigma\sigma^*(t,x)$, it holds that 
$$
\frac{\|a(t,x)\|}{1+|x|^2}+\frac{|b(t,x)|}{1+|x|}\leq C,
$$
where $C>0$ is independent of $t,x$.

(ii) By $({\bf H}^s_{f})$, $({\bf H}^l_{f})$ and $\nu^f_{t,x}(A)=\int_{\mR^d}I_A(f(t,x)z)\nu_1(\dif z)$ for any $A\in\sB(\mR^d)$, it holds that 
\ce
\int_{B_l}|u|^2\nu_{t,x}^f(\dif u)<\infty, ~\mbox{and}~\frac{\int_{B_l}|u|^2\nu_{t,x}^f(\dif u)}{1+|x|^2}\leq C_2,
\de
and
\ce
\nu_{t,x}^f(B^c_l)<\infty, ~\mbox{and}~\int_{B^c_l}\log\left(1+\frac{|u|}{1+|x|}\right)\nu_{t,x}^f(\dif u)\leq C_3.
\de
\er

\br
(i) If $f(t,x)=1$, $({\bf H}^s_{f})$ and $({\bf H}^l_{f})$ become that $\int_{B_l}|z|^2\nu_1(\dif z)<\infty$ and $\nu_1(B^c_l)<\infty$, respectively. These conditions are just right sufficient and necessary for $\nu_1$ to be a L\'evy measure. Therefore, if $f(t,x)\neq1$, it is reasonable to require other conditions.

(ii) If $\nu_1$ is a finite measure, we take $l=0$. And then, Eq.(\ref{Eq1})  goes into 
$$
\dif X_t=b(t,X_t)\dif t+\sigma(t,X_t)\dif B_t+\int_{\mR^d}f(t,X_{t-})zN(\dif t,\dif z).
$$
The type of SDEs has been studied in \cite{q0}. Thus, in the sequel we require $l>0$.
\er

\subsection{Weak solutions of Fokker-Planck equations}

In this subsection, we introduce weak solutions of FPEs.

Consider the FPE associated with Eq.(\ref{Eq1}):
\be
\partial_t\mu_t=\sL_t^*\mu_t,
\label{FPE1}
\ee
where $\sL_t^*$ is the adjoint operator of $\sL_t$, and $(\mu_t)_{t\in[0,T]}$ is a family of probability measures on $\mR^d$.  Weak solutions of Eq.(\ref{FPE1}) are defined as follows. 

\bd\label{weakfpe}
A measurable family $(\mu_t)_{t\in[0,T]}$ of probability measures is called a
 weak solution of the non-local FPE (\ref{FPE1}) if for any $R>0$ and $t\in[0,T]$,
\be
&&\int_0^t\int_{\mR^d}I_{B_R}(x)\left(|b(s,x)|+\|a(s,x)\|+\int_{B_l}|u|^2\nu^f_{s,x}(\dif u)\right)\mu_s(\dif x)\dif s<\infty,
\label{deficond1}\\
&&\int_0^t\int_{\mR^d}\(\nu^f_{s,x}(B_{l\vee(|x|-R)}^c)+I_{B_R}(x)\nu^f_{s,x}(B_l^c)\)\mu_s(\dif x)\dif s<\infty,
\label{deficond2}
\ee
and for all $\phi\in C_c^2(\mR^d)$ and $t\in[0,T]$,
\be
\mu_t(\phi)=\mu_0(\phi)+\int_0^t\mu_s(\sL_s\phi)\dif s,
\label{deficond3}
\ee
where $\mu_t(\phi):=\int_{\mR^d}\phi(x)\mu_t(\dif x)$. The uniqueness of the weak solutions to Eq.(\ref{FPE1}) means that, if $(\mu_t)_{t\in[0,T]}$ and $(\tilde{\mu}_t)_{t\in[0,T]}$ are two weak solutions to Eq.(\ref{FPE1}) with $\mu_0=\tilde{\mu}_0$, then $\mu_t=\tilde{\mu}_t$ for any $t\in[0,T]$.
\ed

By \cite[Remark 1.2]{RXZ}, we know that under the conditions (\ref{deficond1}) (\ref{deficond2}), Eq.(\ref{deficond3}) makes sense. If a weak solution $(\mu_t)_{t\in[0,T]}$ of the non-local FPE (\ref{FPE1}) is absolutely continuous with respect to the Lebesgue measure, then there exists a non-negative measurable function $\rho$ with $\int_{\mR^d}\rho(t,x)\dif x=1$ such that $\mu_t(\dif x)=\rho(t,x)\dif x$. Thus, $\rho$ satisfies the following equation in the distributional sense
\be
\partial_t \rho=-\partial_i(b_i\rho)+\partial_{ij}(a_{ij}\rho)+\sN_t^{f*}\rho.
\label{FPE2}
\ee
Set
$$
\mL:=\left\{\rho\geq 0: \int_{\mR^d}\rho(t,x)\dif x=1, ~\mbox{and}~\sup\limits_{t\in[0,T]}\|\rho(t,\cdot)\|_{\infty}<\infty\right\}.
$$
If there exists a $\rho\in\mL$ satisfying Eq.(\ref{FPE2}) in the distributional sense, we say  Eq.(\ref{FPE2}) has a weak solution in $\mL$.

\subsection{The superposition principle for SDEs driven by L\'evy processes and non-local FPEs}

In the subsection, we state the superposition principle for SDEs driven by L\'evy processes and non-local FPEs. (c.f.\cite[Corollary 1.8]{RXZ})

\bt\label{super}
Suppose that $({\bf H}^1_{b,\sigma})$, $({\bf H}^s_{f})$ and $({\bf H}^l_{f})$ hold and $\mu_0\in\cP(\mR^d)$.

(i) Eq.(\ref{Eq1}) has a martingale solution $\mQ$ with the initial law $\mu_0$ at $s=0$ if and only if Eq.(\ref{FPE1}) has a weak solution $(\mu_t)_{t\in[0,T]}$ starting from $\mu_0$. Moreover, $\mQ\circ e_t^{-1}=\mu_t$ for any $t\in[0,T]$.

(ii) Eq.(\ref{Eq1}) has at most a martingale solution $\mQ$ with the initial law $\mu_0$ at $s=0$ if and only if Eq.(\ref{FPE1}) has at most a weak solution $(\mu_t)_{t\in[0,T]}$ starting from $\mu_0$.
\et

\subsection{The stochastic Gronwall inequality}\label{stogro}
The following stochastic Gronwall inequality comes from \cite[Lemma 2.4]{RXZ}.

\bl\label{stogron}
Let $\xi(t)$ and $\eta(t)$ be two non-negative c\`adl\`ag adapted processes, $A_t$ be a continuous non-decreasing adapted process with $A_0=0$, and $M_t$ be a local martingale with $M_0=0$. Suppose that
\ce
\xi(t)\leq \eta(t)+\int_0^t\xi(s)\dif A_s+M_t, \quad \forall t\geq 0.
\de
Then for any $0<q<p<1$ and any stopping time $\tau>0$, we have
\ce
\left(\mE\left(\sup\limits_{t\in[0,\tau]}\xi(t)^q\right)\right)^{1/q}\leq \left(\frac{p}{p-q}\right)^{1/q}\left(\mE\left(\exp\left\{\frac{pA_{\tau}}{1-p}\right\}\right)\right)^{(1-p)/p}\mE\left(\sup\limits_{t\in[0,\tau]}\eta(t)\right).
\de
\el

\section{The limits of SDEs driven by L\'evy processes}\label{add}

In this section, set $f(t,x):=\gamma g(t,x)$, where $\gamma$ is a real number and $g: [0,T]\times\mR^d\mapsto\mR$ is Borel measurable, and then Eq.(\ref{Eq1}) changes into  
\be
\dif X_t=b(t,X_t)\dif t+\sigma(t,X_t)\dif B_t+\gamma g(t, X_{t-}) \dif L_t, \qquad t\in[0,T].
\label{Eq5}
\ee
Consider the following sequence of SDEs driven by L\'evy processes: for any $n\in\mN$,
\be
\dif X_t^n=b^n(t, X_t^n)\dif t+\sigma^n(t, X_t^n)\dif B_t+\gamma^n g(t, X^n_{t-}) \dif L_t, \quad t\in[0,T],
\label{Eq6}
\ee
where $b^n: [0,T]\times\mR^d\mapsto\mR^d$, $\sigma^n: [0,T]\times\mR^d\mapsto \mR^d\times\mR^m$ are Borel measurable functions and $\{\gamma^n\}$ is a real sequence. When $b^n\rightarrow b, a^n\rightarrow a, \gamma^n\rightarrow\gamma$ in some sense, where $a^n:=\frac{1}{2}\sigma^n\sigma^{n*}$, we study the relationship between martingale solutions of Eq.(\ref{Eq5}) and that of Eq.(\ref{Eq6}). 

The following theorem is the main result in the section.

\bt\label{limit1}
Assume that $b^n, b, \sigma^n, \sigma$ satisfy $({\bf H}^1_{b,\sigma})$ uniformly, $\{\gamma^n\}$ is uniformly bounded, $g$ satisfies $({\bf H}^s_{f})$ and $({\bf H}^l_{f})$, and that Eq.(\ref{FPE2}) has a unique weak solution in $\mL$. Let $\mu_0(\dif x)=\rho_0(x)\dif x\in\cP(\mR^d)$ with $\|\rho_0\|_{\infty}<\infty$, and $\mQ^n, \mQ$ be the martingale solutions of Eq.(\ref{Eq6}) and Eq.(\ref{Eq5}) with the initial law $\mu_0$ at $s=0$, respectively. Assume that

(i) $b^n\rightarrow b, a^n\rightarrow a$ in $L^1_{loc}([0,T]\times\mR^d)$, $\gamma^n\rightarrow\gamma$ as $n\rightarrow\infty$;

(ii) $\mQ^n\circ e_t^{-1}$ is absolutely continuous with respect to the Lebesgue measure on $\mR^d$, $\rho^n(t,x)$ denotes the density, i.e., $\rho^n(t,x):=\frac{(\mQ^n\circ e_t^{-1})(\dif x)}{\dif x}$ for any $t\in[0,T]$  and 
$$
\sup\limits_{t\in[0,T]}\|\rho^n(t,\cdot)\|_{\infty}\leq C,
$$
where $C>0$ is independent of $n$.

Then $\mQ^n\rightarrow\mQ$ in $\cP(D^d_T)$.
\et
\begin{proof}
{\bf Step 1.} We prove that $\{\mQ^n\}_{n\in\mN}$ is tight in $\cP(D^d_T)$.

By Theorem 4.5 in \cite[Page 356]{jjas}, it is sufficient to check that 

(iii) $\lim\limits_{K\rightarrow\infty}\sup\limits_{n}\mQ^n\left(\sup\limits_{t\in[0,T]}|w_t|>K\right)=0$,

(iv) For any stopping time $\tau$, it holds that 
$$
\lim\limits_{\t\rightarrow 0}\sup\limits_{n}\sup\limits_{0\leq\tau<\tau+\t\leq T}\mQ^n\left(|w_{\tau+\t}-w_{\tau}|\geq N\right)=0, \quad \forall N>0.
$$

First of all, \cite[Lemma 3.4]{RXZ} admits us to obtain that there exists a $\psi\in C^2(\mR_+, \mR_+)$ satisfying 
\be
\psi\geq 0,\quad \psi(0)=0, \quad 0<\psi^\prime\leq 1, \quad -2\leq \psi^{\prime\prime}\leq 0, \quad \lim\limits_{r\rightarrow\infty}\psi(r)=\infty,
\label{propsi}
\ee
such that 
\be
\int_{\mR^d}\psi(\log(1+|x|^2))\mu_0(\dif x)<\infty.
\label{inicon}
\ee
Set $\Psi(x):=\psi(\log(1+|x|^2))$, and then $\Psi\in{C^2(\mR^d)}$ with $|\Psi(x)|\leq C\log(2+|x|)$. Note that $\mQ^n$ is a martingale solution of Eq.(\ref{Eq6}) with the initial law $\mu_0$. So, by Remark \ref{equi}, it holds that there exists a local $(\bar{\cB}_t)_{t\in[0,T]}$-adapted martingale $(M^n_t)_{t\in[0,T]}$ under the probability measure $\mQ^n$ such that 
\be
\Psi(w_t)=\Psi(w_0)+\int_0^t\sL^n_s\Psi(w_s)\dif s+M_t^n,
\label{mnlo}
\ee
where $\sL^n_s$ is the infinitesimal generator of Eq.(\ref{Eq6}), i.e.
$$
\sL^n_s\Psi(x)=a^n_{ij}(s,x)\partial_{ij}\Psi(x)+b^n_i(s,x)\partial_i\Psi(x)+\int_{\mR^d}\[\Psi(x+\gamma^n u)-\Psi(x)-I_{|u|\leq l}\gamma^n u^i\partial_i\Psi(x)\]\nu^{g}_{s,x}(\dif u),
$$
and $\nu^{g}_{s,x}(A):=\int_{\mR^d}I_A(g(s,x)z)\nu_1(\dif z)$ for any $A\in\sB(\mR^d)$. 

Next, we estimate $\sL^n_s\Psi(x)$. On one hand, by some calculation, we know that 
$$
\partial_i\Psi(x)=\frac{2x^i}{1+|x|^2}\psi^\prime(\log(1+|x|^2)),
$$
and 
\ce
\partial_{ij}\Psi(x)=\frac{4x^ix^j}{(1+|x|^2)^2}(\psi^{\prime\prime}-\psi^\prime)(\log(1+|x|^2))+\frac{2I_{i=j}}{1+|x|^2}\psi^\prime(\log(1+|x|^2)).
\de
Thus, it follows from Remark \ref{rewrcond} and (\ref{propsi}) that 
\be
a^n_{ij}(s,x)\partial_{ij}\Psi(x)+b^n_i(s,x)\partial_i\Psi(x)\leq C,
\label{esti1}
\ee
where $C>0$ is independent of $n, s, x$. On the other hand, by the mean value theorem, we have that for $|u|\leq l\leq \frac{1}{\sqrt{2}\Gamma}$, where $\Gamma:=\sup\limits_{n}|\gamma^n|$,
\ce
\Psi(x+\gamma^n u)-\Psi(x)-\gamma^n u^i\partial_i\Psi(x)&=&(\gamma^n)^2u^i u^j \partial_{ij}\Psi(x+\delta\gamma^n u)/2\\
&\leq&\Gamma^2\frac{|u|^2}{1+|x+\delta\gamma^n u|^2}\leq \Gamma^2\frac{|u|^2}{1+|x|^2/2- \Gamma^2|u|^2}\\
&\leq& \Gamma^2\frac{2|u|^2}{1+|x|^2},
\de
where $\delta\in[0,1]$, and for $|u|>l$
\ce
\Psi(x+\gamma^n u)-\Psi(x)&=&\psi^\prime(\delta^*)[\log(1+|x+\gamma^n u|^2)-\log(1+|x|^2)]\\
&=&\psi^\prime(\delta^*)\log\left(\frac{1+|x+\gamma^n u|^2}{1+|x|^2}\right)\leq\log\left(1+\frac{|x+\gamma^n u|^2-|x|^2}{1+|x|^2}\right)\\
&=&\log\left(1+\frac{2\gamma^n x^i u^i+|\gamma^n u|^2}{1+|x|^2}\right)\leq\log\left(1+\frac{|\gamma^n u|}{\sqrt{1+|x|^2}}\right)^2\\
&\leq&\log\left(1+\frac{2\Gamma |u|}{1+|x|}\right)^2\leq\log\left(1+\frac{|u|}{1+|x|}\right)^{[4\Gamma]+2}\\
&=&([4\Gamma]+2)\log\left(1+\frac{|u|}{1+|x|}\right),
\de
where $\delta^*\in\mR_+$ and $[4\Gamma]$ stands for the largest integer no more than $4\Gamma$. Thus, it holds that
\ce
&&\int_{\mR^d}\[\Psi(x+\gamma^n u)-\Psi(x)-I_{|u|\leq l}\gamma^n u^i\partial_i\Psi(x)\]\nu^{g}_{s,x}(\dif u)\\
&\leq&2\Gamma^2\frac{\int_{B_l}|u|^2\nu^{g}_{s,x}(\dif u)}{1+|x|^2}+([4\Gamma]+2)\int_{B^c_l}\log\left(1+\frac{|u|}{1+|x|}\right)\nu^{g}_{s,x}(\dif u),
\de
which together with Remark \ref{rewrcond} yields that 
\be
\int_{\mR^d}\[\Psi(x+\gamma^n u)-\Psi(x)-I_{|u|\leq l}\gamma^n u^i\partial_i\Psi(x)\]\nu^{g}_{s,x}(\dif u)\leq C,
\label{esti2}
\ee
where $C>0$ is independent of $n, s, x$. Combining (\ref{esti1}) (\ref{esti2}), we obtain that 
\be
\sL^n_s\Psi(x)\leq C, 
\label{esti3}
\ee
where $C>0$ is independent of $n, s, x$.

Now, inserting (\ref{esti3}) in (\ref{mnlo}), one can have that
\ce
\Psi(w_t)\leq\Psi(w_0)+Ct+M_t^n.
\de
So, Lemma \ref{stogron} admits us to get that
\be
&&\mE^{\mQ^n}\left(\sup\limits_{t\in[0,T]}\Psi^{1/2}(w_t)\right)\leq C\left(\mE^{\mQ^n}\left(\sup\limits_{t\in[0,T]}\(\Psi(w_0)+Ct\)\right)\right)^{1/2}\no\\
&\leq& C\left(\mE^{\mQ^n}\Psi(w_0)+CT\right)^{1/2}=C\left(\int_{\mR^d}\Psi(x)\mu_0(\dif x)+CT\right)^{1/2}\no\\
&=&C\left(\int_{\mR^d}\psi(\log(1+|x|^2))\mu_0(\dif x)+CT\right)^{1/2}<\infty,
\label{esti4}
\ee
where (\ref{inicon}) is used in the last inequality. Thus, (iii) is verified.

For (iv), we have that for any $R>0$, 
\be
&&\mQ^n\left(|w_{\tau+\t}-w_{\tau}|\geq N\right)\no\\
&=&\mQ^n\left(|w_{\tau+\t}-w_{\tau}|\geq N, |w_{\tau}|>R\right)+\mQ^n\left(|w_{\tau+\t}-w_{\tau}|\geq N, |w_{\tau}|\leq R\right)\no\\
&\leq&\mQ^n\left(|w_{\tau}|>R\right)+\mQ^n\left(|w_{\tau+\t}-w_{\tau}|\geq N, |w_{\tau}|\leq R\right)\no\\
&=:&I_1+I_2.
\label{deco}
\ee
For $I_1$, by (\ref{esti4}), it holds that 
\be
I_1&\leq& \frac{\mE^{\mQ^n}\psi^{1/2}(\log(1+|w_\tau|^2))}{\psi^{1/2}(\log(1+R^2))}\leq\frac{\mE^{\mQ^n}\left(\sup\limits_{t\in[0,T]}\Psi^{1/2}(w_t)\right)}{\psi^{1/2}(\log(1+R^2))}\no\\
&\leq& \frac{1}{\psi^{1/2}(\log(1+R^2))}C\left(\int_{\mR^d}\psi(\log(1+|x|^2))\mu_0(\dif x)+CT\right)^{1/2}.
\label{decoi1}
\ee
In the following, we are devoted to dealing with $I_2$. Note that 
\ce
I_2=\mQ^n\left(\mQ^n_{s,y}(|w_{s+\t}-y|\geq N)|_{s=\tau, y=w_{\tau}}, |w_{\tau}|\leq R\right),
\de
where we use the strong Markov property and $\mQ^n_{s,y}:=\mQ^n|_{\cB(D([s,T], \mR^d))}, w(s)=y$ for $s\in[0,T), y\in\mR^d$. Thus, we estimate $\mQ^n_{s,y}(|w_{s+\t}-y|\geq N)$ so as to master $I_2$. 

To treat $\mQ^n_{s,y}(|w_{s+\t}-y|\geq N)$, we define $\Phi(x):=\psi(\log(1+|x-y|^2))$. So, based on Remark \ref{equi}, it holds that there exists a local $(\bar{\cB}_t)_{t\in[0,T]}$-adapted martingale $(\check{M}^n_t)_{t\in[0,T]}$ under the probability measure $\mQ^n_{s,y}$ such that 
\be
\Phi(w_{s+\t})=\Phi(y)+\int_s^{s+\t}\sL^n_r\Phi(w_r)\dif r+\check{M}_{s+\t}^n.
\label{symart}
\ee
Besides, by the similar deduction to that in (\ref{esti1}) (\ref{esti2}), it holds that 
\be
\sL^n_r\Phi(x)&\leq&\frac{2|a^n_{ii}(r,x)|}{1+|x-y|^2}+\frac{2|b^n_i(r,x)||x-y|}{1+|x-y|^2}+2\Gamma^2\frac{\int_{B_l}|u|^2\nu^{g}_{r,x}(\dif u)}{1+|x-y|^2}\no\\
&&+([4\Gamma]+2)\int_{B^c_l}\log\left(1+\frac{|u|}{1+|x-y|}\right)\nu^{g}_{r,x}(\dif u)\no\\
&\leq&C\bigg(\frac{1+|x|^2}{1+|x-y|^2}+\frac{(1+|x|)|x-y|}{1+|x-y|^2}+\frac{1+|x|^2}{1+|x-y|^2}\no\\
&&+\int_{B^c_l}\log\left(1+\frac{|u|}{1+|x-y|}\right)\nu^{g}_{r,x}(\dif u)\bigg)\no\\
&\leq&C(1+|y|^2),
\label{symart1}
\ee  
where we use Remark \ref{rewrcond} in the second inequality. Thus, (\ref{symart}) (\ref{symart1}) yield that 
\be
\Phi(w_{s+\t})\leq\Phi(y)+C(1+|y|^2)\t+\check{M}_{s+\t}^n=C(1+|y|^2)\t+\check{M}_{s+\t}^n.
\label{symart2}
\ee
 Applying Lemma \ref{stogron} to (\ref{symart2}), we have that 
 $$
 \mE^{\mQ^n_{s,y}}\Phi^{1/2}(w_{s+\t})\leq C(1+|y|)\t^{1/2},
 $$
and furthermore
\ce
\mQ^n_{s,y}(|w_{s+\t}-y|\geq N)\leq \frac{\mE^{\mQ^n_{s,y}}\Phi^{1/2}(w_{s+\t})}{\psi^{1/2}(\log(1+N^2))}\leq \frac{C(1+|y|)\t^{1/2}}{\psi^{1/2}(\log(1+N^2))}.
\de
Therefore, it holds that
\be
I_2\leq\frac{C(1+R)\t^{1/2}}{\psi^{1/2}(\log(1+N^2))}.
\label{decoi2}
\ee
Combining (\ref{decoi1}) (\ref{decoi2}), we get that 
\ce
\mQ^n\left(|w_{\tau+\t}-w_{\tau}|\geq N\right)\leq\frac{C\left(\int_{\mR^d}\psi(\log(1+|x|^2))\mu_0(\dif x)+CT\right)^{1/2}}{\psi^{1/2}(\log(1+R^2))}+\frac{C(1+R)\t^{1/2}}{\psi^{1/2}(\log(1+N^2))}.
\de
As $\t\rightarrow0$ and then $R\rightarrow\infty$, one can obtain (iv).

{\bf Step 2.} We show that $\mQ^n$ weakly converges to $\mQ$.

Assume that the limit point of $\{\mQ^n\}_{n\in\mN}$ is $\bar{\mQ}$. And then we only prove that $\mQ=\bar{\mQ}$. Note that $\mQ$ is a martingale solution of Eq.(\ref{Eq5}) with the initial law $\mu_0$, and Eq.(\ref{FPE2}) has a unique weak solution in $\mL$. Thus, by Theorem \ref{super} we further only prove that $\bar{\mQ}$ is a martingale solution of Eq.(\ref{Eq5}) with the initial law $\mu_0$. That is, it is sufficient to check that for $0\leq s<t\leq T$ and a bounded continuous $\bar{\cB}_s$-measurable functional $\chi_s: D^d_T\mapsto \mR$,
\be
\int_{D^d_T}\left[\phi(w_t)-\phi(w_s)-\int_s^t(\sL_r\phi)(w_r)dr\right]\chi_s(w)\bar{\mQ}(\dif w)=0, \quad \forall \phi\in{C_c^2(\mR^d)}.
\label{esti01}
\ee

Next, again note that $\mQ^n\circ e_t^{-1}\rightarrow\bar{\mQ}\circ e_t^{-1}$ in $\cP(\mR^d)$ and $\mQ^n\circ e_0^{-1}=\mu_0=\bar{\mQ}\circ e_0^{-1}$. Thus, by (ii), there exists a $\bar{\rho}(t,x)\geq0$ with $\int_{\mR^d}\bar{\rho}(t,x)\dif x=1$ such that $\bar{\mQ}\circ w_t^{-1}(\dif x)=\bar{\rho}(t,x)\dif x$ and $\rho^n(t,\cdot)\rightarrow\bar{\rho}(t,\cdot)$ in $w^*$-$L^{\infty}(\mR^d)$, where $w^*$-$L^{\infty}(\mR^d)$ is the dual space of $C_c(\mR^d)$, and $\bar{\rho}(0,x)=\rho_0(x)$. Moreover, by the theory of functional analysis and \cite[Lemma 3.8]{RXZ}, we know that for any $\e>0$ and the coefficients $b, a$, there exist $\tilde{b}: [0,T]\times\mR^d\mapsto\mR^d, \tilde{a}: [0,T]\times\mR^d\mapsto \mS_+(\mR^d)$, where $\mS_+(\mR^d)$ is the set of nonnegative definite symmetric $d\times d$ real matrices, and a family of measures $\tilde{\nu}^g_{\cdot,\cdot}$ such that 

(v) $\tilde{b}, \tilde{a}$ are continuous and compactly supported;

(vi) for any $\phi\in{C_c^2(\mR^d)}$, $(t,x)\mapsto\tilde{\sN}^{g1}_t\phi(x)$ and $(t,x)\mapsto\tilde{\sN}^{g2}_t\phi(x)$ are continuous, where $\tilde{\sN}^{g1}_t\phi(x):=\int_{\mR^d}\[\phi(x+\gamma u)-\phi(x)-\gamma\pi^i(u)\partial_i\phi(x)\]\tilde{\nu}^g_{t,x}(\dif u)$, $\tilde{\sN}^{g2}_t\phi(x):=\int_{\mR^d}\[\gamma\pi^i(u)\partial_i\phi(x)-\gamma I_{|u|\leq l} u^i\partial_i\phi(x)\]\tilde{\nu}^g_{t,x}(\dif u)$ and $\pi: \mR^d\mapsto \mR^d$ is a smooth symmetric function with $\pi(u)=u, |u|\leq l$ and $\pi(u)=0, |u|>2l$, and $\sup\limits_{t\in[0,T], x\in\mR^d}|\tilde{\sN}^{g1}_t\phi(x)|<\infty, \sup\limits_{t\in[0,T], x\in\mR^d}|\tilde{\sN}^{g2}_t\phi(x)|<\infty$;

(vii)
\ce
&&\int_0^T\int_{\mR^d}\bigg(|(b_i(r,x)-\tilde{b}_i(r,x))\partial_{i}\phi(x)|+|(a_{ij}(r,x)-\tilde{a}_{ij}(r,x))\partial_{ij}\phi(x)|\\
&&\quad +|\sN^{g1}_r\phi(x)-\tilde{\sN}^{g1}_r\phi(x)|+|\sN^{g2}_r\phi(x)-\tilde{\sN}^{g2}_r\phi(x)|\bigg)\bar{\rho}(r,x)\dif x\dif r<\e.
\de
And then the operator $\tilde{\sL}$ with respect to $\tilde{b}, \tilde{a}, \tilde{\nu}^g_{\cdot,\cdot}$ presents as
\ce
(\tilde{\sL}_r\phi)(x):=\tilde{a}_{ij}(r,x)\partial_{ij}\phi(x)+\tilde{b}_i(r,x)\partial_i\phi(x)+\tilde{\sN}^{g1}_r\phi(x)+\tilde{\sN}^{g2}_r\phi(x), \quad r\in[0,T].
\de

Now, we treat (\ref{esti01}). Note that 
\ce
\int_{D^d_T}\left[\phi(w_t)-\phi(w_s)-\int_s^t(\sL^n_r\phi)(w_r)dr\right]\chi_s(w)\mQ^n(\dif w)=0.
\de
Thus, it holds that
\ce 
&&\left|\int_{D^d_T}\left[\phi(w_t)-\phi(w_s)-\int_s^t(\sL_r\phi)(w_r)dr\right]\chi_s(w)\bar{\mQ}(\dif w)\right|\\
&\leq&\bigg|\int_{D^d_T}\left[\phi(w_t)-\phi(w_s)-\int_s^t(\sL_r\phi)(w_r)dr\right]\chi_s(w)\bar{\mQ}(\dif w)\\
&&\qquad -\int_{D^d_T}\left[\phi(w_t)-\phi(w_s)-\int_s^t(\tilde{\sL}_r\phi)(w_r)dr\right]\chi_s(w)\bar{\mQ}(\dif w)\bigg|\\
&&+\bigg|\int_{D^d_T}\left[\phi(w_t)-\phi(w_s)-\int_s^t(\tilde{\sL}_r\phi)(w_r)dr\right]\chi_s(w)\bar{\mQ}(\dif w)\\
&&\qquad -\int_{D^d_T}\left[\phi(w_t)-\phi(w_s)-\int_s^t(\tilde{\sL}_r\phi)(w_r)dr\right]\chi_s(w)\mQ^n(\dif w)\bigg|\\
&&+\bigg|\int_{D^d_T}\left[\phi(w_t)-\phi(w_s)-\int_s^t(\tilde{\sL}_r\phi)(w_r)dr\right]\chi_s(w)\mQ^n(\dif w)\\
&&\qquad -\int_{D^d_T}\left[\phi(w_t)-\phi(w_s)-\int_s^t(\sL_r\phi)(w_r)dr\right]\chi_s(w)\mQ^n(\dif w)\bigg|\\
&&+\bigg|\int_{D^d_T}\left[\phi(w_t)-\phi(w_s)-\int_s^t(\sL_r\phi)(w_r)dr\right]\chi_s(w)\mQ^n(\dif w)\\
&&\qquad -\int_{D^d_T}\left[\phi(w_t)-\phi(w_s)-\int_s^t(\sL^n_r\phi)(w_r)dr\right]\chi_s(w)\mQ^n(\dif w)\bigg|\\
&=:&J_1+J_2+J_3+J_4.
\de

For $J_1$, we have that
\be
J_1&\leq& C\int_{D^d_T}\left|\int_s^t(\sL_r\phi)(w_r)dr-\int_s^t(\tilde{\sL}_r\phi)(w_r)dr\right|\bar{\mQ}(\dif w)\no\\
&\leq&C\int_s^t\int_{\mR^d}|(\sL_r\phi)(x)-(\tilde{\sL}_r\phi)(x)|\bar{\rho}(r,x)\dif x\dif r\no\\
&\leq&C\int_s^t\int_{\mR^d}\bigg(|(a_{ij}(r,x)-\tilde{a}_{ij}(r,x))\partial_{ij}\phi(x)|+|(b_i(r,x)-\tilde{b}_i(r,x))\partial_{i}\phi(x)|\no\\
&&+|\sN^{g1}_r\phi(x)-\tilde{\sN}^{g1}_r\phi(x)|+|\sN^{g2}_r\phi(x)-\tilde{\sN}^{g2}_r\phi(x)|\bigg)\bar{\rho}(r,x)\dif x\dif r\no\\
&\leq&C\e,
\label{j1es}
\ee
where (vii) is used in the last inequality. For $J_2$, based on the weak convergence of $\{\mQ^n\}$ to $\bar{\mQ}$ and (v) (vi), it holds that there exists a $N_1\in\mN$ such that for $n\geq N_1$ 
\be
J_2\leq \e.
\label{j2es}
\ee
For $J_3$, by the similar deduction to that in $J_1$, one can obtain that 
\ce
J_3&\leq&C\int_s^t\int_{\mR^d}\bigg(|(a_{ij}(r,x)-\tilde{a}_{ij}(r,x))\partial_{ij}\phi(x)|+|(b_i(r,x)-\tilde{b}_i(r,x))\partial_{i}\phi(x)|\no\\
&&+|\sN^{g1}_r\phi(x)-\tilde{\sN}^{g1}_r\phi(x)|+|\sN^{g2}_r\phi(x)-\tilde{\sN}^{g2}_r\phi(x)|\bigg)\rho^n(r,x)\dif x\dif r.
\de
So, Remark \ref{rewrcond} and (ii) (v), together with the Fatou lemma, yield that there exists a $N_2\in\mN, N_2\geq N_1$ such that for $n\geq N_2$
\be
J_3\leq C\e.
\label{j3es}
\ee
For $J_4$, we get that
\ce
J_4&\leq& C\int_{D^d_T}\left|\int_s^t(\sL_r\phi)(w_r)dr-\int_s^t(\sL^n_r\phi)(w_r)dr\right|\mQ^n(\dif w)\no\\
&\leq&C\int_s^t\int_{\mR^d}|(\sL_r\phi)(x)-(\sL^n_r\phi)(x)|\rho^n(r,x)\dif x\dif r\no\\
&\leq&C\int_s^t\int_{\mR^d}\bigg(|(a_{ij}(r,x)-a^n_{ij}(r,x))\partial_{ij}\phi(x)|+|(b_i(r,x)-b^n_i(r,x))\partial_{i}\phi(x)|\bigg)\rho^n(r,x)\dif x\dif r\no\\
&&+C\int_s^t\int_{\mR^d}\bigg|\int_{\mR^d}\[\phi(x+\gamma u)-\phi(x)-\gamma I_{|u|\leq l} u^i\partial_i\phi(x)\]\nu^g_{t,x}(\dif u)\no\\
&&\qquad\qquad\qquad -\int_{\mR^d}\[\phi(x+\gamma^n u)-\phi(x)-\gamma^n I_{|u|\leq l} u^i\partial_i\phi(x)\]\nu^g_{t,x}(\dif u)\bigg|\rho^n(r,x)\dif x\dif r,
\de
and furthermore by (i) and the Fatou lemma, there exists a $N_3\geq N_2$ such that for $n\geq N_3$
\be
J_4\leq C\e.
\label{j4es}
\ee
Combining (\ref{j1es})-(\ref{j4es}), one can obtain that
\ce
\left|\int_{D^d_T}\left[\phi(w_t)-\phi(w_s)-\int_s^t(\sL_r\phi)(w_r)dr\right]\chi_s(w)\bar{\mQ}(\dif w)\right|\leq C\e.
\de
Letting $\e\rightarrow 0$, we have (\ref{esti01}). The proof is complete.
\end{proof}

\section{The robustness of the nonlinear filterings}\label{robunon}

In this section, set $g(t,x)=1$ in Eq.(\ref{Eq5})-(\ref{Eq6}) and then Eq.(\ref{Eq5})-(\ref{Eq6}) change into  
\be
\dif X_t=b(t,X_t)\dif t+\sigma(t,X_t)\dif B_t+\gamma\dif L_t, \qquad t\in[0,T],
\label{Eq7}
\ee
and 
\be
\dif X_t^n=b^n(t, X_t^n)\dif t+\sigma^n(t, X_t^n)\dif B_t+\gamma^n\dif L_t, \quad t\in[0,T],
\label{Eq8}
\ee
respectively. We define nonlinear filtering problems associated with Eq.(\ref{Eq7}) and Eq.(\ref{Eq8}) and then study the relationship of two nonlinear filterings under the framework of Theorem \ref{limit1}.

\subsection{Nonlinear filtering problems}\label{nonfil}

In the subsection, we introduce nonlinear filtering problems associated with Eq.(\ref{Eq7}) and Eq.(\ref{Eq8}).

Given the filtered probability space $(\Omega, \mathscr{F}, \{\mathscr{F}_t\}_{t\in[0,T]},\mP)$. Let $B_{\cdot}, L_{\cdot}$ be $m$-dimensional Brownian motion and $d$-dimensional pure jump L\'evy process defined on it, respectively. We assume:

\begin{enumerate}[({\bf H}$^{2}_{b,\sigma}$)]
\item There exist two positive constants $C_b, C_\sigma$ such that
for any $t\in[0,T]$ and $x,y\in\mR^{d}$
\ce
&&|b(t,x)-b(t,y)|\leq C_b|x-y|\cdot
\log(|x-y|^{-1}+\me);\\
&&\|\sigma(t,x)-\sigma(t,y)\|^2\leq C_\sigma|x-y|^{2}\cdot
\log(|x-y|^{-1}+\me).
\de
\end{enumerate}

Under the assumption ({\bf H}$^{1}_{b,\sigma}$)-({\bf H}$^2_{b,\sigma}$), it holds that Eq.(\ref{Eq7}) and Eq.(\ref{Eq8}) have unique strong solutions denoted as $(X_t)$  and $(X^n_t)$ with $\mP\circ X_0^{-1}=\mu_0$ and $\mP\circ (X^n_0)^{-1}=\mu_0$, respectively. 

In the following, we introduce the nonlinear filtering problem associated with $(X_t)$. Given an observation process $(Y_t)_{t\in[0,T]}$ as follows:
\ce
Y_t=W_t+\int_0^t h(X_s)\dif s+\int_0^t\int_{\mU_0}u\tilde{N}_{\lambda}(\dif s, \dif u)+\int_0^t\int_{\mR^k\setminus\mU_0}u N_{\lambda}(\dif s, \dif u),
\de
where $W_{\cdot}$ is a $k$-dimensional Brownian motion and $N_{\lambda}(\dif t,\dif u)$ is a random measure with a predictable compensator $\lambda(X_t,u)\dif t\nu_2(\dif u)$. Here the function $\lambda: \mR^d\times\mR^k\mapsto(0,1)$ is Borel measurable and $\nu_2$ is a $\sigma$-finite measure defined on $\mR^k$ with $\nu_2(\mR^k\setminus\mU_0)<\infty$ and $\int_{\mU_0}|u|^2\,\nu_2(\dif u)<\infty$ for a fixed $\mU_0\in\sB(\mR^k)$. Concretely speaking, set 
$$
\tilde{N}_\lambda(\dif t, \dif u):=N_\lambda(\dif t, \dif u)-\lambda(X_t,u)\dif t\nu_2(\dif u), \quad t\in[0,T], 
$$ 
and then $\tilde{N}_\lambda(\dif t, \dif u)$ is the compensated martingale measure
of $N_{\lambda}(\dif t, \dif u)$. Moreover, we require that $B_{\cdot}, L_{\cdot}, W_{\cdot}, N_{\lambda}(\dif t, \dif u)$ are mutually independent. $h: \mR^d\mapsto\mR^k$ is Borel measurable. Here,
we assume more:

\begin{enumerate}[(\bf{H}$^1_h$)]
\item 
\ce
\int_0^T|h(X_t)|^2\dif t<\infty, ~\mbox{and}~ \int_0^T|h(X^n_t)|^2\dif t<\infty.
\de
\end{enumerate}
\begin{enumerate}[(\bf{H}$^1_{\lambda}$)]
\item There exists a positive function $L(u)$ satisfying
\ce
\int_{\mU_0}\frac{\left(1-L(u)\right)^2}{L(u)}\nu_2(\dif u)<\infty,
\de
such that $0<\iota \leq L(u)<\lambda(x,u)<1$ for $u\in\mU_0$, where $\iota$ is a constant.
\end{enumerate}

Now, denote
\ce
\Sigma_t^{-1}:&=&\exp\bigg\{-\int_0^t h^i(X_s)\dif W^i_s-\frac{1}{2}\int_0^t
\left|h(X_s)\right|^2\dif s-\int_0^t\int_{\mU_0}\log\lambda(X_{s-},u)\tilde{N}_{\lambda}(\dif s, \dif u)\\
&&\quad\qquad -\int_0^t\int_{\mU_0}\(1-\lambda(X_s,u)+\lambda(X_s,u)\log\lambda(X_s,u)\)\nu_2(\dif u)\dif s\bigg\}.
\de
Thus, by ({\bf H}$^1_h$) ({\bf H}$^1_{\lambda}$) we know that $\Sigma_{\cdot}^{-1}$ is an exponential martingale. Define a measure $\tilde{\mP}$ via
$$
\frac{\dif \tilde{\mP}}{\dif \mP}=\Sigma_T^{-1}.
$$
Under the probability measure $\tilde{\mP}$, it follows from the Girsanov theorem that $\tilde{W}_{\cdot}:=W_{\cdot}+\int_0^{\cdot} h(X_s)\dif s$ is a Brownian motion and 
$$
\eta_{\cdot}:=\int_0^{\cdot}\int_{\mU_0}u\tilde{N}_{\lambda}(\dif s, \dif u)+\int_0^{\cdot}\int_{\mR^k\setminus\mU_0}u N_{\lambda}(\dif s, \dif u)
$$
is a pure jump L\'evy process with the L\'evy measure $\nu_2$. Moreover, $X_{\cdot}$ is independent of $\tilde{W}_{\cdot}, \eta_{\cdot}$ under the probability measure $\tilde{\mP}$. And then we rewrite $\Sigma_t$ as
\ce
\Sigma_t&=&\exp\bigg\{\int_0^t h^i(X_s)\dif \tilde{W}^i_s-\frac{1}{2}\int_0^t
\left|h(X_s)\right|^2\dif s+\int_0^t\int_{\mU_0}\log\lambda(X_{s-},u)N_{\lambda}(\dif s, \dif u)\\
&&\quad\qquad +\int_0^t\int_{\mU_0}\(1-\lambda(X_s,u)\)\nu_2(\dif u)\dif s\bigg\}.
\de
Set
\ce
&&\varrho_t(\phi):=\mE^{\tilde{\mP}}[\phi(X_t)\Sigma_t|\mathscr{F}_t^{Y}],\\
&&\pi_t(\phi):=\mE[\phi(X_t)|\mathscr{F}_t^{Y}], \qquad \phi\in\cB(\mR^d),
\de
where $\mE^{\tilde{\mP}}$ stands for the expectation under the probability measure $\tilde{\mP}$ and $\mathscr{F}_t^{Y} \triangleq\sigma(Y_s:
 0\leq s \leq t)$. And then by the Kallianpur-Striebel formula it holds that
\ce
\pi_t(\phi)=\frac{\varrho_t(\phi)}{\varrho_t(1)}.
\de

Next, we introduce the nonlinear filtering problem associated with $(X^n_t)$. Set
\ce
Y^n_t:=W_t+\int_0^t h(X^n_s)\dif s+\int_0^t\int_{\mU_0}u\tilde{N}^n_{\lambda}(\dif s, \dif u)+\int_0^t\int_{\mR^k\setminus\mU_0}u N^n_{\lambda}(\dif s, \dif u),
\de
where $N^n_{\lambda}(\dif s,\dif u)$ is a random measure with a predictable compensator $\lambda(X^n_s,u)\dif s\nu_2(\dif u)$ and $\tilde{N}^n_\lambda(\dif s, \dif u):=N^n_\lambda(\dif s, \dif u)-\lambda(X^n_s,u)\dif s\nu_2(\dif u)$. 
By the similar way to above, we can define $\Sigma_t^{n}, \tilde{\mP}^n, \tilde{W}^n, \eta^n, \varrho^n_t$ and $\pi^n_t$ by replacing $X_t, Y_t$ with $X^n_t, Y^n_t$. Moreover, we require that $B_{\cdot}, L_{\cdot}, W_{\cdot}, N^n_{\lambda}(\dif t, \dif u)$ are mutually independent. Here, we remind that $\pi^n_t, \pi_t$ are defined under the same probability measure $\mP$.

\subsection{The relationship between $\pi^n_t$ and $\pi_t$}\label{relpi}

In the subsection, we observe the relationship between $\pi^n_t$ and $\pi_t$ under the framework of Theorem \ref{limit1}. 

First of all, under the assumption of Theorem \ref{limit1} and ({\bf H}$^2_{b,\sigma}$), by Theorem \ref{limit1} and \cite[Theorem 2.26, Page 157]{jjas} we know that $\mP\circ X_{\cdot}^{-1}=\mQ$, $\mP\circ (X^n_{\cdot})^{-1}=\mQ^n$ and $X_{\cdot}^n\Longrightarrow X_{\cdot}$ in $\cP(D_T^d)$, where ``$\Longrightarrow$" denotes convergence in distribution of random variables as well as weak convergence of probability measures. And then we apply some functionals to prove that $\pi_{\cdot}^n\Longrightarrow\pi_{\cdot}$ in $\cP(D([0,T],\cP(\mR^d)))$. To do this, we assume more:
\begin{enumerate}[(\bf{H}$^2_h$)]
\item $h$ is continuous and satisfies
\ce
\lim\limits_{n\rightarrow\infty}\mE\left(\int_0^T|h(X^n_s)-h(X_s)|^2\dif s\right)=0
\de
\end{enumerate}
\begin{enumerate}[(\bf{H}$^2_{\lambda}$)]
\item $\lambda$ is continuous in the first variable $x$ and satisfies
\ce
\lim\limits_{n\rightarrow\infty}\mE\left(\int_0^T\int_{\mU_0}|\log\lambda(X^n_s, u)-\log\lambda(X_s, u)|^2\nu_2(\dif u)\dif s\right)=0.
\de
\end{enumerate}

Thus, by the assumptions ({\bf H}$^2_h$) ({\bf H}$^2_{\lambda}$), one can obtain that 
\be
\left(X_{\cdot}^n, Z^n_{\cdot}, V^n_{\cdot}\right)\Longrightarrow\left(X_{\cdot}, Z_{\cdot}, V_{\cdot}\right),
\label{indiscon}
\ee
where
\ce
&&Z^n_t:=\int_0^t |h(X^n_s)|^2\dif s, \quad Z_t:=\int_0^t |h(X_s)|^2\dif s,\\
&&V^n_t:=\int_0^t\int_{\mU_0}\(1-\lambda(X^n_s,u)+\log\lambda(X^n_s,u)\)\nu_2(\dif u)\dif s,\\
&&V_t:=\int_0^t\int_{\mU_0}\(1-\lambda(X_s,u)+\log\lambda(X_s,u)\)\nu_2(\dif u)\dif s.
\de
Next, note that by the Skorokhod representation theorem, there exist a probability space $(\Omega^0, \sF^0, \mP^0)$ and $\bar{X}_{\cdot}^n, \bar{Z}^n_{\cdot}, \bar{V}^n_{\cdot}, \bar{X}_{\cdot}, \bar{Z}_{\cdot}, \bar{V}_{\cdot}$ on it such that
\be
\left(\bar{X}_{\cdot}^n, \bar{Z}^n_{\cdot}, \bar{V}^n_{\cdot}\right)\rightarrow\left(\bar{X}_{\cdot}, \bar{Z}_{\cdot}, \bar{V}_{\cdot}\right)  \quad a.s. \mP^0,
\label{ascon}
\ee
and 
\ce
\cL\left(\bar{X}_{\cdot}^n, \bar{Z}^n_{\cdot}, \bar{V}^n_{\cdot}\right)=\cL\left(X_{\cdot}^n, Z^n_{\cdot}, V^n_{\cdot}\right), \quad \cL\left(\bar{X}_{\cdot}, \bar{Z}_{\cdot}, \bar{V}_{\cdot}\right)=\cL\left(X_{\cdot}, Z_{\cdot}, V_{\cdot}\right),
\de
where $\cL$ denotes the joint distribution. Besides, let $\Omega^1:=C([0,T],\mR^k)$, $\sF^1$ be the Borel $\sigma$-field on $\Omega^1$  and $\mP^1$ be the Wiener measure on $(\Omega^1, \sF^1)$. Let $\bar{W}$ be the canonical process on $(\Omega^1, \sF^1, \mP^1)$. Let $\Omega^2:=D([0,T],\mR^k)$. And then we equip $\Omega^2$ with the Skorokhod topology and  $\sF^2$ denotes the Borel $\sigma$-field induced by the Skorokhod topology. Moreover, we take $\mP^2=\tilde{\mP}\circ \eta_{\cdot}^{-1}$ and then $(\Omega^2, \sF^2, \mP^2)$ is a probability space. $\bar{\eta}$ denotes the canonical process on it. Let $(\bar{\Omega}, \bar{\sF}, \bar{\mP}):=(\Omega^0, \sF^0, \mP^0)\times(\Omega^1, \sF^1, \mP^1)\times(\Omega^2, \sF^2, \mP^2)$. We remind that the distribution of $(X_{\cdot}, \tilde{W}_{\cdot}, \eta_{\cdot})$ on $(\Omega, \sF, \tilde{\mP})$ is the same to that of $(\bar{X}_{\cdot}, \bar{W}_{\cdot}, \bar{\eta}_{\cdot})$ on $(\bar{\Omega}, \bar{\sF}, \bar{\mP})$, and the distribution of $(X^n_{\cdot}, \tilde{W}^n_{\cdot}, \eta^n_{\cdot})$ on $(\Omega, \sF, \tilde{\mP}^n)$ is the same to that of $(\bar{X}_{\cdot}^n, \bar{W}_{\cdot}, \bar{\eta}_{\cdot})$ on $(\bar{\Omega}, \bar{\sF}, \bar{\mP})$.

In the following, we present $\pi^n_{\cdot}, \pi_{\cdot}$ as some functionals on $(\bar{\Omega}, \bar{\sF}, \bar{\mP})$. Set 
\ce
&&\<F_t(w^1, w^2), \phi\>:=\int_{\Omega^0}\phi(\bar{X}_t(w^0))q_t(w^0, w^1, w^2)\mP^0(\dif w^0), \\
&&\<F^n_t(w^1, w^2), \phi\>:=\int_{\Omega^0}\phi(\bar{X}^n_t(w^0))q^n_t(w^0, w^1, w^2)\mP^0(\dif w^0), \quad \phi\in\cB(\mR^d),
\de
where 
\ce
q_t(w^0, w^1, w^2)&:=&\exp\bigg\{\int_0^t h^i(\bar{X}_s(w^0))\dif \bar{W}^i_s-\frac{1}{2}\int_0^t
\left|h(\bar{X}_s(w^0))\right|^2\dif s\\
&&\quad\qquad +\int_0^t\int_{\mU_0}\log\lambda(\bar{X}_{s-}(w^0),u)\tilde{N}_{\kappa}(\dif s, \dif u)\\
&&\quad\qquad +\int_0^t\int_{\mU_0}\(1-\lambda(\bar{X}_s(w^0),u)+\log\lambda(\bar{X}_{s}(w^0),u)\)\nu_2(\dif u)\dif s\bigg\},\\
q^n_t(w^0, w^1, w^2)&:=&\exp\bigg\{\int_0^t h^i(\bar{X}^n_s(w^0))\dif \bar{W}^i_s-\frac{1}{2}\int_0^t
\left|h(\bar{X}^n_s(w^0))\right|^2\dif s\\
&&\quad\qquad +\int_0^t\int_{\mU_0}\log\lambda(\bar{X}^n_{s-}(w^0),u)\tilde{N}_{\kappa}(\dif s, \dif u)\\
&&\quad\qquad +\int_0^t\int_{\mU_0}\(1-\lambda(\bar{X}^n_s(w^0),u)+\log\lambda(\bar{X}^n_{s}(w^0),u)\)\nu_2(\dif u)\dif s\bigg\},
\de
and 
$$
\kappa_t:=w^2_t-w^2_{t-}, \quad N_{\kappa}((0, t], A):=\#\{0<s\leq t, \kappa_s\in A\}, \quad A\in\sB(\mR^k\setminus\{0\}),
$$
and $\tilde{N}_{\kappa}(\dif t, \dif u):=N_{\kappa}(\dif t, \dif u)-\nu_2(\dif u)\dif t$ is the compensated martingale measure of the Poisson random measure $N_{\kappa}(\dif t, \dif u)$. And then it holds that
\ce
\<F_t(\tilde{W}_{\cdot}, \eta_{\cdot}), \phi\>=\varrho_t(\phi), \quad \<F^n_t(\tilde{W}_{\cdot}^n, \eta_{\cdot}^n), \phi\>=\varrho^n_t(\phi).
\de
Moreover, by the similar deduction to that on the top of Theorem 3.2 in \cite{qd}, one can get that 
$$
\<F_t(w^1, w^2), 1\>>0, \quad \<F^n_t(w^1, w^2), 1\>>0, \quad a.s.\bar{\mP}.
$$
Thus, we define 
\ce
\<H_t(w^1, w^2), \phi\>:=\frac{\<F_t(w^1, w^2), \phi\>}{\<F_t(w^1, w^2), 1\>}, \quad \<H^n_t(w^1, w^2), \phi\>:=\frac{\<F^n_t(w^1, w^2), \phi\>}{\<F^n_t(w^1, w^2), 1\>}
\de
and then obtain 
\ce
\<H_t(\tilde{W}_{\cdot}, \eta_{\cdot}), \phi\>=\pi_t(\phi), \quad \<H^n_t(\tilde{W}_{\cdot}^n, \eta_{\cdot}^n), \phi\>=\pi^n_t(\phi).
\de

Now, it is the position to state and prove the main result in the section.
\bt\label{filcon}
Under the assumption of Theorem \ref{limit1} and ({\bf H}$^2_{b,\sigma}$) ({\bf H}$^1_{h}$)-({\bf H}$^2_{h}$) ({\bf H}$^1_{\lambda}$)-({\bf H}$^2_{\lambda}$), it holds that $\mP\circ(\pi^n)_{\cdot}^{-1}\Longrightarrow\mP\circ\pi_{\cdot}^{-1}$ in $\cP(D([0,T], \cP(\mR^d)))$.
\et
\begin{proof}
First of all, note that for $G\in C_b(D([0,T], \cP(\mR^d)))$, 
\ce
&&\mE[G(\pi_{\cdot}^n)]=\mE^{\tilde{\mP}^n}\left[G(H_{\cdot}^n(\tilde{W}_{\cdot}^n, \eta_{\cdot}^n))\lambda^n_T\right]=\mE^{\bar{\mP}}\left[G(H_{\cdot}^n)q_T^n\right],\\
&&\mE[G(\pi_{\cdot})]=\mE^{\tilde{\mP}}\left[G(H_{\cdot}(\tilde{W}_{\cdot}, \eta_{\cdot}))\lambda_T\right]=\mE^{\bar{\mP}}\left[G(H_{\cdot})q_T\right].
\de
Therefore, we only need to prove that as $n\rightarrow\infty$, $H_{\cdot}^n\rightarrow H_{\cdot}$ and $q_T^n\rightarrow q_T$ in the probability measure $\bar{\mP}$.

Next, we are devoted to showing $H_{\cdot}^n\rightarrow H_{\cdot}$ in the probability measure $\bar{\mP}$. And then by the definition of $H_{\cdot}^n, H_{\cdot}$, it is sufficient to prove that $F_{\cdot}^n\rightarrow F_{\cdot}$ in the probability measure $\bar{\mP}$. This is implied by for any $t_n\rightarrow t$ and any $\e>0$,
\be
\lim\limits_{n\rightarrow\infty}(\mP^1\times\mP^2)\{|\<F_{t_n}^n(w^1, w^2), \phi\>-\<F_t(w^1, w^2), \phi\>|\geq \e\}=0, \quad \forall \phi\in C_b(\mR^d).
\label{guijie}
\ee
Note that 
\be
&&(\mP^1\times\mP^2)\left\{|\<F_{t_n}^n(w^1, w^2), \phi\>-\<F_t(w^1, w^2), \phi\>|\geq \e\right\}\no\\
&\leq&(\mP^1\times\mP^2)\left\{\int_{\Omega^0}|\phi(\bar{X}^n_{t_n}(w^0))q^n_{t_n}(w^0, w^1, w^2)-\phi(\bar{X}_t(w^0))q_t(w^0, w^1, w^2)|\mP^0(\dif w^0)\geq \e\right\}\no\\
&\leq&(\mP^1\times\mP^2)\left\{\int_{\Omega^0}|\phi(\bar{X}^n_{t_n}(w^0))q^n_{t_n}(w^0, w^1, w^2)-\phi(\bar{X}^n_{t_n}(w^0))q^n_t(w^0, w^1, w^2)|\mP^0(\dif w^0)\geq \e/4\right\}\no\\
&&+(\mP^1\times\mP^2)\left\{\int_{\Omega^0}|\phi(\bar{X}^n_{t_n}(w^0))q^n_t(w^0, w^1, w^2)-\phi(\bar{X}^n_{t_n}(w^0))q_t(w^0, w^1, w^2)|\mP^0(\dif w^0)\geq \e/4\right\}\no\\
&&+(\mP^1\times\mP^2)\left\{\int_{\Omega^0}|\phi(\bar{X}^n_{t_n}(w^0))q_t(w^0, w^1, w^2)-\phi(\bar{X}^n_t(w^0))q_t(w^0, w^1, w^2)|\mP^0(\dif w^0)\geq \e/4\right\}\no\\
&&+(\mP^1\times\mP^2)\left\{\int_{\Omega^0}|\phi(\bar{X}^n_t(w^0))q_t(w^0, w^1, w^2)-\phi(\bar{X}_t(w^0))q_t(w^0, w^1, w^2)|\mP^0(\dif w^0)\geq \e/4\right\}\no\\
&=:&\Sigma_1+\Sigma_2+\Sigma_3+\Sigma_4.
\label{estimate}
\ee
So, we estimate $\Sigma_1, \Sigma_2, \Sigma_3, \Sigma_4$ to obtain (\ref{guijie}).

For $\Sigma_1$, by Lemma \ref{expro} below and the dominated convergence theorem, we know that 
\ce
\lim\limits_{n\rightarrow\infty}\int_{\bar{\Omega}}|\phi(\bar{X}^n_{t_n}(w^0))q^n_{t_n}(w^0, w^1, w^2)-\phi(\bar{X}^n_{t_n}(w^0))q^n_t(w^0, w^1, w^2)|\dif \bar{\mP}(w^0, w^1, w^2)=0.
\de
By Chebychev's inequality, it holds that for any $\delta>0$, there exists a $N_1\in\mN$ such that for $n>N_1$, 
\be
\Sigma_1\leq \delta/4. \label{estimate1}
\ee

For $\Sigma_2$, by ({\bf H}$^2_h$) ({\bf H}$^2_{\lambda}$) and (\ref{ascon}), it holds that 
\ce
&&\int_0^t h^i(\bar{X}^n_s(w^0))\dif \bar{W}^i_s\rightarrow \int_0^t h^i(\bar{X}_s(w^0))\dif \bar{W}^i_s, \quad ~\mbox{in}~\bar{\mP},\\
&&\int_0^t\left|h(\bar{X}^n_s(w^0))\right|^2\dif s\rightarrow \int_0^t\left|h(\bar{X}_s(w^0))\right|^2\dif s, \quad a.s.\bar{\mP},\\
&&\int_0^t\int_{\mU_0}\log\lambda(\bar{X}^n_{s-}(w^0),u)\tilde{N}_{\kappa}(\dif s, \dif u) \rightarrow \int_0^t\int_{\mU_0}\log\lambda(\bar{X}_{s-}(w^0),u)\tilde{N}_{\kappa}(\dif s, \dif u), \quad ~\mbox{in}~\bar{\mP},\\
&&\int_0^t\int_{\mU_0}\(1-\lambda(\bar{X}^n_s(w^0),u)+\log\lambda(\bar{X}^n_{s}(w^0),u)\)\nu_2(\dif u)\dif s\\
&&\quad\rightarrow\int_0^t\int_{\mU_0}\(1-\lambda(\bar{X}_s(w^0),u)+\log\lambda(\bar{X}_{s}(w^0),u)\)\nu_2(\dif u)\dif s, a.s.\bar{\mP}.
\de
Thus, we know that $q^n_t(w^0, w^1, w^2)\rightarrow q_t(w^0, w^1, w^2)$ in the probability measure $\bar{\mP}$, which together with $\int_{\bar{\Omega}}q^n_t(w^0, w^1, w^2) \dif \bar{\mP}(w^0, w^1, w^2)=1, \int_{\bar{\Omega}}q_t(w^0, w^1, w^2) \dif \bar{\mP}(w^0, w^1, w^2)=1$ and the Scheffe Lemma, yields that 
\ce
\lim\limits_{n\rightarrow\infty}\int_{\bar{\Omega}}|q^n_t(w^0, w^1, w^2)-q_t(w^0, w^1, w^2)|\dif \bar{\mP}(w^0, w^1, w^2)=0,
\de
and furthermore
\ce
\lim\limits_{n\rightarrow\infty}\int_{\bar{\Omega}}|\phi(\bar{X}^n_{t_n}(w^0))q^n_t(w^0, w^1, w^2)-\phi(\bar{X}^n_{t_n}(w^0))q_t(w^0, w^1, w^2)|\dif \bar{\mP}(w^0, w^1, w^2)=0.
\de
From this, it follows that there exists a $N_2\in\mN, N_2\geq N_1$ such that for $n>N_2$, 
\be
\Sigma_2 \leq \delta/4. \label{estimate2}
\ee

For $\Sigma_3$, note that $\bar{X}^n_{t_n}(w^0)\rightarrow \bar{X}^n_t(w^0)$ in the probability measure $\bar{\mP}$ (\cite[Definition 1.6, Page 3]{sa}). Thus, by the dominated convergence theorem, we have that 
\ce
\lim\limits_{n\rightarrow\infty}\int_{\bar{\Omega}}|\phi(\bar{X}^n_{t_n}(w^0))q_t(w^0, w^1, w^2)-\phi(\bar{X}^n_t(w^0))q_t(w^0, w^1, w^2)|\dif \bar{\mP}(w^0, w^1, w^2)=0.
\de
And then it follows from Chebychev's inequality that there exists a $N_3\in\mN, N_3\geq N_2$ such that for $n>N_3$, 
\be
\Sigma_3\leq \delta/4.\label{estimate3}
\ee

For $\Sigma_4$, by (\ref{ascon}) and the dominated convergence theorem, it holds that 
\ce
\lim\limits_{n\rightarrow\infty}\int_{\bar{\Omega}}|\phi(\bar{X}^n_t(w^0))q_t(w^0, w^1, w^2)-\phi(\bar{X}_t(w^0))q_t(w^0, w^1, w^2)|\dif \bar{\mP}(w^0, w^1, w^2)=0.
\de
So, we get that there exists a $N_4\in\mN, N_4\geq N_3$ such that for $n>N_4$, 
\be
\Sigma_4\leq \delta/4.
\label{estimate4}
\ee

Combining (\ref{estimate1})-(\ref{estimate4}) with (\ref{estimate}), one can obtain that for $n>N_4$, 
\ce
(\mP^1\times\mP^2)\left\{|\<F_{t_n}^n(w^1, w^2), \phi\>-\<F_t(w^1, w^2), \phi\>|\geq \e\right\}\leq \delta.
\de
Thus, (\ref{guijie}) is proved. 

Finally, by the similar deduction to that about $\Sigma_2$, we have that $q_T^n\rightarrow q_T$ in the probability measure $\bar{\mP}$. So, the proof is complete.
\end{proof}

\bl\label{expro}
$\int_0^t\int_{\mU_0}\log\lambda(\bar{X}_{s-}(w^0),u)\tilde{N}_{\kappa}(\dif s, \dif u)$ is stochastic continuous in $t$.
\el
\begin{proof}
Note that 
\ce
&&\int_0^T\int_{\mU_0}|\log\lambda(\bar{X}_{s}(w^0),u)|^2\nu_2(\dif u)\dif s
\leq\int_0^T\int_{\mU_0}|\log L(u)|^2\nu_2(\dif u)\dif s\\
&\leq& \int_0^T\int_{\mU_0}\frac{(1-L(u))^2}{L^2(u)}\nu_2(\dif u)\dif s
\leq\int_0^T\int_{\mU_0}\frac{(1-L(u))^2}{L(u)}\frac{1}{\iota}\nu_2(\dif u)\dif s<\infty.
\de
Thus, by \cite[Theorem 4.2.12, Page 228]{da}, we know that $\int_0^t\int_{\mU_0}\log\lambda(\bar{X}_{s-}(w^0),u)\tilde{N}_{\kappa}(\dif s, \dif u)$ is right continuous in $t$ and then right stochastic continuous in $t$.

Besides, we take $t_n\uparrow t$ as $n\rightarrow\infty$ for $t_n, t\in[0,T]$. And then 
\ce
\lim\limits_{n\rightarrow\infty}\int_{t_n}^t\int_{\mU_0}|\log\lambda(\bar{X}_{s}(w^0),u)|^2\nu_2(\dif u)\dif s=0, a.s. \bar{\mP}.
\de
So, for any $\delta, \eta>0$, there exists a $N\in\mN$ such that for $n>N$,
\ce
\bar{\mP}\left\{\int_{t_n}^t\int_{\mU_0}|\log\lambda(\bar{X}_{s}(w^0),u)|^2\nu_2(\dif u)\dif s>\delta\right\}<\eta,
\de
which together with \cite[Exercise 4.2.10, Page 228]{da}, yields that for any $\e>0$ such that
\ce
&&\bar{\mP}\left\{\left|\int_0^{t_n}\int_{\mU_0}\log\lambda(\bar{X}_{s-}(w^0),u)\tilde{N}_{\kappa}(\dif s, \dif u)-\int_0^t\int_{\mU_0}\log\lambda(\bar{X}_{s-}(w^0),u)\tilde{N}_{\kappa}(\dif s, \dif u)\right|\geq \e\right\}\\
&\leq&\frac{\delta}{\e^2}+\bar{\mP}\left\{\int_{t_n}^t\int_{\mU_0}|\log\lambda(\bar{X}_{s}(w^0),u)|^2\nu_2(\dif u)\dif s>\delta\right\}\\
&<&\frac{\delta}{\e^2}+\eta.
\de
From this, it follows that $\int_0^t\int_{\mU_0}\log\lambda(\bar{X}_{s-}(w^0),u)\tilde{N}_{\kappa}(\dif s, \dif u)$ is left stochastic continuous in $t$. The proof is complete.
\end{proof}

\section{The appendix}\label{app}

{\bf Verification of Remark \ref{equi}.}

{\bf Necessity.} First of all, we choose a smooth function $\chi_n$ such that $\chi_n(x)=1, |x|\leq n$ and $\chi_n(x)=0, |x|\geq 2n$. And then for any $\phi\in C^2(\mR^d)$ with $|\phi(x)|\leq C\log(2+|x|)$, $\phi_n:=\phi\chi_n\in C_c^2(\mR^d)$. From this, it follows that 
\ce
\cM^{\phi_n}_t=\phi_n(w_t)-\phi_n(w_s)-\int_s^t(\sL_r\phi_n)(w_r)\dif r
\de
is a $(\bar{\cB}_t)_{t\in[0,T]}$-adapted martingale under the probability measure $\mQ$. Set $\tau_v=\inf\{T\geq t\geq s, |w_t|>v\}$ for $v\in\mN$, and then $\{\tau_v\}$ is a $(\bar{\cB}_t)_{t\in[0,T]}$-stopping time sequence and $\tau_v\uparrow T$ as $v\rightarrow\infty$. Thus, 
\ce
\cM^{\phi_n}_{t\land\tau_v}&=&\phi_n(w_{t\land\tau_v})-\phi_n(w_{s\land\tau_v})-\int_{s\land\tau_v}^{t\land\tau_v}(\sL_r\phi_n)(w_r)\dif r
\de
is still a $(\bar{\cB}_t)_{t\in[0,T]}$-adapted martingale under $\mQ$. The dominated convergence theorem admits us to obtain
\ce
\cM^{\phi}_{t\land\tau_v}&=&\phi(w_{t\land\tau_v})-\phi(w_{s\land\tau_v})-\int_{s\land\tau_v}^{t\land\tau_v}(\sL_r\phi)(w_r)\dif r
\de
is also a $(\bar{\cB}_t)_{t\in[0,T]}$-adapted martingale under $\mQ$. That is, 
\ce
\cM^{\phi}_t=\phi(w_t)-\phi(w_s)-\int_s^t(\sL_r\phi)(w_r)\dif r
\de
is a $(\bar{\cB}_t)_{t\in[0,T]}$-adapted local martingale under $\mQ$.

{\bf Sufficiency.}  For any $\phi\in C_c^2(\mR^d)$, we know that
\ce
\cM^{\phi}_t=\phi(w_t)-\phi(w_s)-\int_s^t(\sL_r\phi)(w_r)\dif r
\de
is a $(\bar{\cB}_t)_{t\in[0,T]}$-adapted local martingale under the probability measure $\mQ$. So, there exists a $(\bar{\cB}_t)_{t\in[0,T]}$-stopping time sequence $\{\tau_n, n\in\mN\}$ with $\tau_n\uparrow T$ such that 
\ce
\cM^{\phi}_{t\land\tau_n}&=&\phi(w_{t\land\tau_n})-\phi(w_{s\land\tau_n})-\int_{s\land\tau_n}^{t\land\tau_n}(\sL_r\phi)(w_r)\dif r
\de
is a $(\bar{\cB}_t)_{t\in[0,T]}$-adapted martingale under $\mQ$. By the dominated convergence theorem, it holds that
\ce
\cM^{\phi}_t=\phi(w_t)-\phi(w_s)-\int_s^t(\sL_r\phi)(w_r)\dif r
\de
is a $(\bar{\cB}_t)_{t\in[0,T]}$-adapted martingale under $\mQ$. The proof is complete.

\bigskip

\textbf{Acknowledgements:}

The author is very grateful to Professor Xicheng Zhang and Renming Song for valuable discussions.

\end{document}